\documentclass{jnmp01b}

\usepackage{amsmath}

\numberwithin{equation}{section}

\def\dps{\displaystyle}
\def \uline{\underbar}

\def\co{{\cal {O}}}
\def\pmq{\pmb{$q$}}
\def\pmtheta{\pmb{\ensuremath{\theta}}}
\def\N{{\bf {N}}}
\def\Z{{\bf {Z}}}

\newenvironment{remark}[1]%
{\begin{trivlist}\item[]\textbf{Remark #1.}\ }
{\end{trivlist}}

\makeatletter
\DeclareRobustCommand{\primfrac}[1]{%
  \PackageWarning{amsmath}{%
Foreign command \@backslashchar#1; %
\protect\frac\space or \protect\genfrac\space should be used instead%
  }
  \global\@xp\let\csname#1\@xp\endcsname\csname @@#1\endcsname
  \csname#1\endcsname
}
\makeatother

\begin{document}

\renewcommand{\evenhead}{B.A. Kupershmidt}
\renewcommand{\oddhead}{\pmq-Newton Binomial: \ From  Euler To Gauss}


\thispagestyle{empty}

\begin{flushleft}
\footnotesize \sf
Journal of Nonlinear Mathematical Physics \qquad 2000, V.7, N~2,
\pageref{firstpage}--\pageref{lastpage}.
\hfill {\sc (Re)-Examining the Past}
\end{flushleft}

\vspace{-5mm}

\copyrightnote{2000}{B.A. Kupershmidt}

\Name{\pmq-Newton Binomial: \ From  Euler To Gauss}

\label{firstpage}

\Author{Boris A. KUPERSHMIDT}

\Adress{The University of Tennessee Space Institute,
Tullahoma, TN  37388, USA \\
E-mail: \ bkupersh@utsi.edu}

\Date{Received March 6, 2000; Accepted April 7, 2000}

\begin{abstract}
\noindent
A counter-intuitive result of Gauss (formulae (1.6), (1.7)
below) is made less mysterious by virtue of being generalized through the
introduction of an additional parameter.
\end{abstract}


\section{A formula of Gauss revisited}

Consider the Newton binomial for a positive integer $N$:
\begin{equation*}
(1-x)^N = \sum^N_{\ell =0} {N \choose \ell} (1 - x)^\ell. \tag{1.1}
\end{equation*}
Substituting $x=1$ into this formula, we get
\begin{equation*}
\sum^N_{\ell =0} {N \choose \ell} (-1)^\ell = 0 . \tag{1.2}
\end{equation*}

What happens with these two equalities in the $q$-mathematics framework?  Newton's formula (1) 
becomes Euler's formula
\begin{equation*}
(1 -^{\hskip-.09truein \cdot}x  )^N = (1 - x) (1 - qx) ... 
(1 - q^{N-1} x) = \sum^N_{\ell = 0} \bigg[{N \atop 
\ell}\bigg] (- x)^\ell q^{{\ell \choose 2}} , \tag{1.3}
\end{equation*}
where $\dps{\bigg[{N \atop \ell}\bigg] = \bigg[{N \atop \ell}\bigg]_q}$ are the Gaussian 
polynomials, or $q$-binomial coefficients:
\begin{align*}
\bigg[&{n \atop k}\bigg] = {[n]! \over [k]! [n-k]!} = {[n]...[n-k+1] \over [k]!} , 
\ \ \ k \in \Z_+, \tag{1.4}
\\[1ex]
[&k]! = [k]_q! = [1] [2] ... [k], \ \ \ [0]! = 1, \ \ \ [n] = [n]_q = (1 - q^n) (1 - q). 
\end{align*}
Substituting $x=1$ into the Euler formula (1.3), we find
\begin{equation*}
\sum^N_{\ell = 0} \bigg[{N \atop \ell}\bigg] (- 1)^\ell q^{{\ell \choose 2}} = 0. 
\tag{1.5}
\end{equation*}
This does not look exactly as a $q$-analogue of formula (1.2).  \\
How about the sum $\dps{\sum^N_{\ell = 0} 
\bigg[{N \atop \ell}\bigg] (-1)^\ell }$? 

The answer is quite surprising.  Denote
\begin{equation*}
s_{N|0} = (-1)^N \sum^N_{\ell =0} \bigg[{N \atop \ell}\bigg] (-1)^\ell. \tag{1.6}
\end{equation*}
Gauss found that
\begin{gather*}
s_{2m+1|0} = 0, \ \ \ m \in \Z_+, \tag{1.7a}
\\
s_{2m+2|0} = (1-q) (1 - q^3) ... (1-q^{2m+1}), \ \ m \in \Z_+. \tag{1.7b}
\end{gather*}
These formulae are easy to prove, but they are nevertheless mystifying:  there is no hint in the 
definition (1.6) that some sort of 
2-periodicity is involved.  In addition, formula (1.2) may claim the 
following sums as proper $q$-analogues:
\begin{equation*}
s_{N|1} = (-1)^N \sum^N_{\ell =0} \bigg[{N \atop \ell}\bigg](-q)^\ell, \tag{1.8}
\end{equation*}
or even
\begin{equation*}
s_{N|r} = (-1)^N \sum^N_{\ell =0} \bigg[{N \atop \ell}\bigg] (-q^r)^\ell. \tag{1.9}
\end{equation*}
Indeed, we shall verify later on that
\begin{gather*}
s_{2m+1|1} = - (1 - q^{2m+1}) s_{2m|0} = - \sqcap^m_{t=0} (1-q^{2t+1}), \tag{1.10a}
\\
s_{2m|1} = s_{2m|0} = \sqcap^m_{t=1} (1 - q^{2t-1}). \tag{1.10b}
\end{gather*}
Similar but more complex formulae can be derived for other values of $r$, not just for $r=0$ 
and $r=1$.  We shall abstain from such derivations, as they are superseded by the general 
formulae (1.12) below.

What seems to be happening here is that the functions
\begin{equation*}
S_N (x) = (-1)^N \sum^N_{\ell =0} \bigg[{N \atop \ell}\bigg] (-x)^\ell \tag{1.11}
\end{equation*}
possess some interesting properties worthy of attention; and once the decision to pay attention has 
been made, one quickly conjectures the formulae

\begin{align*}
S_{2m+1} (x) & = - \sum^{2m+1}_{\ell =0} \bigg[{2m+1 \atop \ell} \bigg] 
(-x)^\ell \\ 
& = \sum^m_{k=0} \bigg[ {m \atop k}\bigg]_{q^{2}} (x -^{\hskip-.09truein \cdot}
  1)^{2m+1-2k} (q^{2m+1} ; 
q^{-2})_k, \tag{1.12a}  
\end{align*}
\begin{align*}
S_{2m+2} (x) & = \sum^{2m+2}_{\ell =0} \bigg[{2m+2 \atop \ell}\bigg] 
(-x)^\ell \\
& = \sum^{m+1}_{k=0} \bigg[{m+1 \atop k}\bigg]_{q{^2}} (x -^{\hskip-.09truein \cdot} 
1)^{2m+2-2k} (q^{2m+1} ; 
q^{-2})_k. \tag{1.12b} 
\end{align*}
The additional notations employed above are to be understood as
\begin{equation*}
(u \dot + v)^\ell = \sqcap^{\ell-1}_{k=0} (u + q^k v), \ \ \ell \in \N; \ \ \ (u \dot + 
v)^0 = 1, \tag{1.13}
\end{equation*}
and
\begin{equation*}
(a; Q)_\ell = \sqcap^{\ell-1}_{k=0} (1 - Q^k a), \ \ \ \ell \in \N; \ \ \  (a ; Q)_0 
= 1. \tag{1.14}
\end{equation*}
If we define
\begin{equation*}
\epsilon (N) = \bigg\{ { \ 1, \ \ N \ {\rm{is \ even}} \atop 0, \ \ N \ {\rm{is \ odd}} }
= \bigg\lfloor{N+2 \over 2}\bigg\rfloor - \bigg\lfloor{N+1 \over 2}\bigg\rfloor, 
\tag{1.15}
\end{equation*}
then formulae (1.12) can be rewritten as
\begin{align*}
S_N(x) & = (-1)^N \sum^N_{\ell =  0} \bigg[{N \atop \ell} \bigg] (-x)^\ell \\
& = \sum^{\lfloor N/2 \rfloor}_{k=0} \bigg[{\lfloor N/2 \rfloor \atop k} 
\bigg]_{q^{2}} (x -^{\hskip-.09truein \cdot} 1)^{N-2k} (q^{N-\epsilon 
(N)}; q^{-2})_k. \tag{1.16} 
\end{align*}

Substituting $x=1$ into formulae (1.12) we recover Gauss' formulae (1.7).

Let us now prove formulae (1.12).  Denote the RHS of formulae (1.16) by $\tilde S_N(x)$.  To 
show that
\begin{equation*}
S_N (x) = \tilde S_N (x), \tag{1.17}
\end{equation*}
we shall verify, first, that
\begin{gather*}
{dS_N \over d_qx} = [N] S_{N-1}, \tag{1.18a}
\\
{d \tilde S_N \over d_q x} = [N] \tilde S_{N-1} (x), \tag{1.18b}
\end{gather*}
and second, that
\begin{equation*}
S_N (1) = \tilde S_N (1); \tag{1.19}
\end{equation*}
here 
\begin{equation*}
{d f(x) \over d_q x}  = {f (qx) - f(x) \over qx - x} \tag{1.20}
\end{equation*}
is the $q$-derivative.  Since $S_1 = \tilde S_1 = x -1$, these verifications would suffice.

We start with formula (1.18a).  We have:
\begin{align*}
{dS_N \over d_q x} 
&= {d \over d_q x} \bigg((-1)^N \sum^N_{\ell =0} \bigg[{N \atop \ell}\bigg] 
(-1)^\ell x^\ell \bigg) 
= (-1)^N \sum^N_{\ell =1} \bigg[{N \atop \ell}\bigg] [\ell] 
(-1)^\ell x^{\ell-1} \ \ [{\rm{by}} \ (1.22)]  
\\
&=(-1)^N [N] \sum^N_{\ell=1} \bigg[{N-1 \atop \ell -1} \bigg] 
(-x)^{\ell -1} (-1) 
= [N] (-1)^{N-1} \sum^{N-1}_{\ell=0} \bigg[{N-1 \atop \ell}\bigg] 
(-x)^\ell 
\\
&= [N] S_{N-1}, \tag{1.21}
\end{align*}
where we used the obvious formula
\begin{equation*}
\bigg[{w \atop \ell} \bigg] [\ell] = [w] \bigg[{w-1 \atop \ell -1}\bigg]. \tag{1.22}
\end{equation*}

Next, formula (1.18b), which we shall check separately for odd and even $N$, making use of the 
easy verifiable relation 
\begin{equation*}
{d (x \dot + v)^\alpha \over d_q x} = [\alpha] (x \dot + v)^{\alpha -1}. \tag{1.23}
\end{equation*}
So, for $N$ odd, we have
\begin{equation*}
{d \tilde S_1 \over d_qx} = {d \over d_q x} (x-1) = 1 = \tilde S_0, \tag{1.24}
\end{equation*}
and then
\begin{align*}
{d \tilde S_{2m+3} \over d_qx} &= \sum^{m+1}_{k=0} \bigg[{m+1 \atop k}\bigg]_{q^{2}} [2m+3 - 
2k] (x -^{\hskip-.09truein \cdot}   1)^{2m+2 -2k} (q^{2m+3}; q^{-2} )_k 
\\
&= [2m+3] \sum^{m+1}_{k=0} \bigg[{m+1 \atop k}\bigg]_{q^{2}} (x -^{\hskip-.09truein \cdot}
  1)^{2m+2 - 2k} (q^{2m+1} 
; q^{-2})_k = [2m+3] \tilde S_{2m+2}, 
\end{align*}
because
\begin{equation*}
[2m+3 -2k] (q^{2m+3}; q^{-2})_k = [2m+3] (q^{2m+1}; q^{-2})_k \ ; \tag{1.25}
\end{equation*}
for $N$ even, we find
\begin{align*}
{d \tilde S_{2m+2} \over d_qx} &= \sum^{m+1}_{k=0} \bigg[{m +1 \atop k}\bigg]_{q^{2}} [2m+2 - 
2k] (x \dot -1)^{2m+1-2k} (q^{2m+1}; q^{-2})_k 
\\
&= [2m+2] \sum^m_{k=0} \bigg[{m \atop k}\bigg]_{q^{2}} (x -^{\hskip-.09truein \cdot}
 1)^{2m+1-2k} (q^{2m+1}; 
q^{-2})_k = [2m+2] \tilde S_{2m+1}, 
\end{align*}
because
\begin{equation*}
\bigg[{m+1 \atop k}\bigg]_{q^{2}} [2m+2-2k] = [2m+2] \bigg[{m \atop k}\bigg]_{q^{2}}, 
\tag{1.26}
\end{equation*}
which is true in view of the obvious relation
\begin{equation*}
[u]_{q^{2}} = [2u]_q/[2]_q. \tag{1.27}
\end{equation*}

It remains to verify formula (1.19), which is nothing but the Gauss formula (1.7).  We shall 
verify the latter in 4 easy steps. \\
\uline{$1^{st}$ Step}  is formula (1.7a):
\begin{align*}
s_{2m+1|0} & = - \sum_{\ell \geq 0} \bigg[{2m+1 \atop \ell}\bigg]
 (-1)^\ell = \sum_{\ell \geq 0} 
\bigg[{2m+1 \atop 2m+1-\ell}\bigg] (-1)^{\ell-1} \\
& =  \sum_{L \geq 0} \bigg[{2m+1 \atop L}\bigg] (-1)^{2m-L} =  - s_{2m+1|0}, 
\end{align*}
so that $s_{2m+1|0}= 0; $\\
\uline{$2^{nd}$ Step}  is formula (1.10b):
\begin{equation*}
s_{2m|1} = \sum_{\ell \geq 0} \bigg[{2m \atop \ell}\bigg] (-q)^{\ell} = \sum_{\ell \geq 0} 
\bigg[{2m \atop \ell}\bigg] (-1)^{\ell} = s_{2m|0}. \tag{1.28}
\end{equation*}
Indeed, 
\begin{equation*}
{s_{2m|1} - s_{2m|0} \over q-1} = \sum_{\ell \geq 0} \bigg[{2m \atop \ell } \bigg] 
(-1)^{\ell} [\ell] \ [{\rm{by}} \ (1.22)] = - [2m] \sum_{\ell \geq 1} \bigg[{2m-1 \atop \ell -1}
\bigg] (-1)^{\ell -1} 
\end{equation*}
${\rm{[by \ (1.7a)]}} = 0; $ \\
\uline{$3^{rd}$ Step}  is formula (1.10a): 
\begin{equation*}
\sum_{\ell \geq 0} \bigg[{2m+1 \atop \ell}\bigg] (-q)^{\ell} = (1 - q^{2m+1} )
 \sum_{\ell \geq 
0} \bigg[{2m \atop \ell}\bigg] (-1)^\ell. \tag{1.29}
\end{equation*}
Indeed, since
\begin{equation*}
\bigg[{2m+1 \atop \ell}\bigg] = \bigg[{2m \atop \ell} \bigg] + q^{2m+1-\ell} \bigg[{2m 
\atop \ell -1}\bigg], \tag{1.30}
\end{equation*}
we have:
\begin{gather*}
\sum_{\ell \geq 0} \bigg[{2m+1 \atop \ell}\bigg] (-q)^\ell = \sum_{\ell \geq  0} (-q)^\ell 
\bigg[{2m \atop \ell}\bigg] + \sum_{\ell \geq 1} (-q)^\ell q^{2m+1-\ell} \bigg[{2m \atop 
\ell -1} \bigg] \ {\rm{[by}} \ (1.28) ] 
\\
\qquad= s_{2m|0} - q^{2m+1} \sum \bigg[{2m \atop \ell -1}\bigg] (-1)^{\ell -1} = (1 - q^{2m+1}) 
s_{2m|0} ; 
\end{gather*} 
\uline{$4^{th}$ Step} is the last one:  we prove that
\begin{equation*}
s_{2m+2|0} = (1 - q^{2m+1} ) s_{2m|0}, \tag{1.31}
\end{equation*}
from which the Gauss formula (1.76) follows at once, since
\begin{equation*}
s_{2|0} = 1 - [2] + 1 = 1 - (1 + q) +1 = 1-q. \tag{1.32}
\end{equation*}
Now, 
\begin{align*}
s_{2m+2|0} &= \sum_{\ell \geq 0} \bigg[{2m+2 \atop \ell}\bigg] (-1)^{\ell} \ {\rm{[by}} \ 
(1.28)] = \sum_{\ell \geq 0} \bigg[{2m+2 \atop \ell}\bigg] (-q)^{\ell} \ \ {\rm{[by}} \ 
(1.30)] 
\\
&= \sum_{\ell \geq 0} (-q)^{\ell} \bigg[{2m+1 \atop \ell}\bigg] + \sum_{\ell \geq 1} 
(-q)^{\ell} q^{2m+2-\ell} \bigg[{2m+1 \atop \ell -1}\bigg] \ {\rm{[by}} \ (1.29) ] 
\\
&= (1-q^{2m+1}) s_{2m|0} - q^{2m+2} \sum_{\ell \geq 1} \bigg[{2m+1 \atop \ell -1} \bigg] 
(-1)^{\ell -1} \ \ {\rm{[by \ (1.7a)]}} \\
&= (1-q^{2m+1})s_{2m|0}. 
\end{align*}

We are done.  Formula (1.16) is thereby proven.  Substituting into this formula $x=0$, we get 
an interesting identity 
\begin{equation*}
\sum^{\lfloor N/2 \rfloor}_{k=0} \bigg[{\lfloor N/2 \rfloor \atop k }\bigg]
_{q^{2}} q^{({N-2k \atop 2})}
 (q^{N-\epsilon (N)} ; 
q^{-2} )_k = 1. \tag{1.33}
\end{equation*}

\section{A different proof}

To prove {\it{polynomial}} identities (1.12) generalizing Gauss' formulae (1.7), we had to 
prove independently the Gauss result along the way.  This is not entirely agreeable.  One 
ought to prove formulae (1.12) directly, by-passing the verification of the original Gauss 
formulae.

Such a proof follows.

Let $R_N(x)$ stand for either $S_N(x)$ or $\tilde S_N(x)$.  We shall verify that
\begin{equation*}
R_{N+1} (x) = x R_N (x) - R_N (qx). \tag{2.1}
\end{equation*}
Since
\begin{equation*}
S_0 (x) = \tilde S_0 (x) = 1, \ \ \ S_1 (x) = \tilde S_1 (x) = x - 1, \tag{2.2}
\end{equation*}
such a verification will prove that $S_N (x) = \tilde S_N (x) $ for all $N$.

We start with $R_N (x) = S_N (x)$.  Let's look for a relation of the form
\begin{equation*}
S_{N+1} (x) = Bx S_N (bx) + AS_N (ax) . \tag{2.3}
\end{equation*}
The $x^0$ - coefficients (recall that $S_n (x) = (-1)^n \sum_{\ell} \bigg[\dps{
{n \atop \ell}}\bigg] 
(-x)^\ell)$ yield
\begin{equation*}
A = -1; \tag{2.4}
\end{equation*}
The $x^{N+1} $- coefficients yield
\begin{equation*}
Bb^N = 1 \ \Leftrightarrow \ B = b^{-N}; \tag{2.5}
\end{equation*}
Finally, for $0 < r < N + 1$, the $x^r$ - coefficients provide
\begin{equation*}
\bigg[{N+1 \atop r} \bigg] = Bb^{r-1} \bigg[{N \atop r-1}\bigg] + a^r \bigg[{N \atop r}\bigg] 
. \tag{2.6}
\end{equation*}
In view of the relation (2.5), formula (2.6) can be rewritten as
\begin{equation*}
\bigg[{N+1 \atop r}\bigg] = (b^{-1})^{N+1-r} \bigg[{N \atop r-1}\bigg] + a^r \bigg[{N \atop 
r}\bigg]. \tag{2.7}
\end{equation*}
Now, since
\begin{align*}
\bigg[{N+1 \atop r}\bigg] & = \bigg[{N \atop r-1}\bigg] 
+ q^r \bigg[{N \atop r}\bigg]  \tag{2.8a}  \\
&= q^{N+1-r} \bigg[{N \atop r-1}\bigg] + \bigg[{N \atop r}\bigg], \tag{2.8b} 
\end{align*}
equation (2.7) has two solutions:
\begin{gather*}
b=1, \ \ a=q, \tag{2.9a}
\\
b=q^{-1}, \ \ a=1. \tag{2.9b}
\end{gather*}
Thus,
\begin{align*}
S_{N+1} (x) & = x S_N (x) - S_N (qx)   \tag{2.10a} \\
 & = q^N x S_N (q^{-1} x) - S_N (x).   \tag{2.10b} 
\end{align*}
(For $q=1$, we get just {\it{one}} relation, $S_{N+1} (x) = (x - 1) S_N (x).) $

Denote by $\co $ the linear operator acting on functions of $x$ by the rule:
\begin{equation*}
\co (f(x)) = x f(x) - f (qx). \tag{2.11}
\end{equation*}
We need to check that
\begin{equation*}
\co (\tilde S_N) = \tilde S_{N+1}. \tag{2.12}
\end{equation*}
We shall check separately the cases of even and odd $N$:
\begin{gather*}
\tilde S_{2m+1} (x) = \sum^m_{k=0} (x -^{\hskip-.09truein \cdot}
1)^{2m+1-2k} c_{m|k}, \tag{2.13}
\\
c_{m|k} = \bigg[{m \atop k}\bigg]_{q^{2}} (q^{2m+1}; q^{-2})_k, \tag{2.14}
\\
\tilde S_{2m} (x) = \sum^m_{k=0} (x -^{\hskip-.09truein \cdot}
 1)^{2m-2k} d_{m|k}, \tag{2.15}
\\
d_{m|k} = \bigg[{m \atop k}\bigg]_{q^{2}} (q^{2m-1}; q^{-2})_k. \tag{2.16}
\end{gather*}

To proceed further, let's establish first that
\begin{equation*}
\co ((x -^{\hskip-.09truein \cdot} 
1)^s) = (x -^{\hskip-.09truein \cdot} 1)^{s+1} + q^{s-1} (1 - q^s) (x \dot -1 )^{s-1} . 
\tag{2.17}
\end{equation*}
Indeed, 
\begin{align*}
\co ((x \dot -1 )^s) &= x (x -^{\hskip-.09truein \cdot} 
1)^s - (q x -^{\hskip-.09truein \cdot}  1)^s = ((x - q^s) + q^s) (x \dot -1)^s 
- q^s (x -^{\hskip-.09truein \cdot}  1 q^{-1})^s  
\\
&= (x-q^s) (x -^{\hskip-.09truein \cdot} 1)^s + q^s (x -^{\hskip-.09truein \cdot} 1
)^{s-1}  (x - q^{s-1}) - q^s (x-q^{-1}) 
(x -^{\hskip-.09truein \cdot}  1) 
^{s-1}  
\\
&= (x -^{\hskip-.09truein \cdot} 1)^{s+1} + q^s (x \dot -1 )^{s-1} 
((x - q^{s-1}) - (x - q^{-1})) \\
&= (x -^{\hskip-.09truein \cdot}  1) 
^{s+1} + q^s (x -^{\hskip-.09truein \cdot} 1)^{s-1} q^{-1} (1 - q^s). 
\end{align*}

Now, 
\begin{align*}
\co (\tilde S_{2m+1}) &= \sum^m_{k=0} c_{m|k} \co((x -^{\hskip-.09truein \cdot}
  1) ^{2m+1-2k})  
\\
&= \sum^m_{k=0} c_{m|k} ((x -^{\hskip-.09truein \cdot}  1)^{2m+2-2k} + q^{2m-2k} (1 - q^{2m+1-2k})
 (x -^{\hskip-.09truein \cdot}x   1)
^{2m - 2k})  
\\
&= \sum^{m+1}_{k=0} (x -^{\hskip-.09truein \cdot}  1)^{2m+2 - 2k} 
(c_{m|k} + c_{m|k-1} q^{2m+2 -2k} (1-q^{2m+3-2k})), \tag{2.18$\ell$}
\end{align*}
while
\begin{equation*}
\tilde S_{2m+2} = \sum^{m+1}_{k=0} (x -^{\hskip-.09truein \cdot} 1)^{2m+2 - 2k} d_{m+1|k}, \tag{2.18r}
\end{equation*}
so we need to verify that
\begin{equation*}
d_{m+1|k} = c_{m|k} + c_{m|k-1} q^{2m+2-2k} (1 - q^{2m+3-2k}), \tag{2.19}
\end{equation*}
which is
\begin{gather*}
\qquad= \bigg[{m+1 \atop k}\bigg]_{q^{2}} (q^{2m+1}; q^{-2})_k 
\\
\qquad=\bigg[{m \atop k}\bigg]_{q^{2}} 
(q^{2m+1}; q^{-2})_k + \bigg[{m \atop k-1}\bigg]_{q^{2}} (q^{2m+1}; q^{-2})_{k-1} q^{
2m+2-2k} (1-q^{2m+3-2k}), \tag{2.20}
\end{gather*}
which is equivalent to
\begin{equation*}
\bigg[{m + 1 \atop k} \bigg]_{q^{2}} = \bigg[{m \atop k}\bigg]_{q^{2}} + \bigg[{m \atop 
k-1}\bigg]_{q^{2}} (1 - q^{2m+1-2(k-1)} )^{-1} q^{2m+2-2k} (1-q^{2m+3-2k}), 
\notag
\end{equation*}
which is finally
\begin{equation*}
\bigg[{m+1 \atop k}\bigg]_{q^{2}} = \bigg[{m \atop k}\bigg]_{q^{2}} 
+ \bigg[{m \atop k -1}\bigg]_{q^{2}} (q^2)^{m+1-k}, 
\end{equation*}
and this is so by formula (2.8b).

Next,
\begin{align*}
\co (\tilde S_{2m}) &= \sum^m_{k=0} d_{m|k} \co ((x -^{\hskip-.09truein \cdot}
 1)^{2m-2k}) 
\\
&= \sum^m_{k=0} d_{m|k} ((x -^{\hskip-.09truein \cdot} 1)^{2m+1-2k} + q^{2m-2k-1} (1 - q^{2m-2k}) 
(x -^{\hskip-.09truein \cdot}  1)^{
2m-1-2k} ) 
\\
&= \sum^m_{k=0} (x -^{\hskip-.09truein \cdot}
 1)^{2m+1-2k} (d_{m|k} + d_{m|k-1} q^{2m+1-2k} (1-q^{2m+2-2k})), 
\tag{2.21$\ell$}
\end{align*}
while
\begin{equation*}
\tilde S_{2m+1} = \sum^m_{k=0} (x -^{\hskip-.09truein \cdot}  1)^{2m+1-2k} c_{m|k}, \tag{2.21r}
\end{equation*}
so we need to check that
\begin{equation*}
c_{m|k} = d_{m|k} + d_{m|k-1} q^{2m+1-2k} (1 - q^{2m+2-2k}), \tag{2.22}
\end{equation*}
which is
\begin{gather*}
\bigg[{m \atop k}\bigg]_{q^{2}} (q^{2m+1}; q^{-2})_k = \bigg[{m \atop k}\bigg]_{q^{2}} 
(q^{2m-1}; q^{-2} )_k  
\\
\qquad+\bigg[{m \atop k-1}\bigg]_{q^{2}} (q^{2m-1}; q^{-2})_{k-1} 
q^{2m+1-2k} (1-q^{2m+2-2k}), \tag{2.23}
\end{gather*}
which is equivalent to
\begin{equation*}
\bigg[{m \atop k}\bigg]_{q^{2}}(1 - q^{2m+1}) = \bigg[{m \atop k}\bigg]_{q^{2}} 
(1-q^{2m-1-2(k-1)}) + \bigg[{m \atop k-1}\bigg]_{q^{2}} q^{2m+1-2k} (1-q^{2m+2-2k} ), 
\notag
\end{equation*}
which can be rewritten as
\begin{equation*}
\bigg[{m \atop k}\bigg]_{q^{2}}
 q^{2m+1-2k} (1-q^{2k}) = \bigg[{m \atop k-1}\bigg]_{q^{2}} q^{2m+1-2k} 
(1-q^{2m+2-2k} ), 
\end{equation*}
which is equivalent to
\begin{equation*}
\bigg[{m \atop k}\bigg] [k] = \bigg[{m \atop k-1}\bigg] [m +1-k], 
\end{equation*}
which is obvious.
\begin{remark}{2.24}
Set
\begin{equation*}
\tilde S_N (x) = \sum_k e_{N|k} (x -^{\hskip-.09truein \cdot}1)^{N -2k}, \tag{2.25}
\end{equation*}
so that
\begin{equation*}
c_{m|k} = e_{2m+1|k}, \ \ d_{m|k} = e_{2m|k}. \tag{2.26}
\end{equation*}
Then the pair of equalities (2.19) and (2.22) can be rewritten as the single one:
\begin{equation*}
e_{N+1|k} = e_{N|k} + e_{N|k-1} q^{N+1-2k} (1 - q^{N+2-2k}), \tag{2.27}
\end{equation*}
equivalent to the relation
\begin{equation*}
\tilde S_{N+1} = \co (\tilde S_N).
\end{equation*}
\end{remark}

\section{The Taylor expansions point of view}

Formula (1.16) (or (2.25)) is reminiscent of the Taylor expansion:
\begin{equation*}
f(x) = \sum_{k \geq 0} {f^{(k)} (a) \over k!} (x -a)^k, \tag{3.1}
\end{equation*}
where
\begin{equation*}
f^{(i)} (x) = \bigg({d \over dx}\bigg)^i (f (x)). \tag{3.2}
\end{equation*}
There exist many different $q$-versions of the classical Taylor expansion.  We shall make 
use below of the following particular one:
\begin{equation*}
f(x) = \sum_{k \geq 0} {f^{(k)} (a) \over [k]!} (x -^{\hskip-.09truein \cdot}
   a)^k, \tag{3.3}
\end{equation*}
where now
\begin{equation*}
f^{(k)} (x) = \bigg({d \over d_q x}\bigg)^k (f(x)). \tag{3.4}
\end{equation*}

We shall prove formula (3.3) for $f$ being polynomial in $x$.  It's enough to consider the
case $f(x) = x^n$, so that
\begin{equation*}
f^{(k)} (x) = [k]! \bigg[{n \atop k}\bigg] x^{n-k}, \tag{3.5}
\end{equation*}
and we thus have to check that
\begin{equation*}
x^n = \sum_k \bigg[{n \atop k}\bigg] a^{n-k} (x -^{\hskip-.09truein \cdot}
   a)^k. \tag{3.6}
\end{equation*}
This can be verified either directly, or deduced from the identity 
(formula (2.10) in [5], p.~75)
\begin{equation*}
\sum_k \bigg[{n \atop k}\bigg] a^{n-k} (x \dot + b)^k = \sum_k \bigg[{n \atop k} \bigg] 
x^{n-k} (a \dot + b)^k. \tag{3.7}
\end{equation*}
for $b = -a$.

Thus, formula (3.3) is proven.  Taking $f(x)$ to be $S_N (x)$, 
\begin{equation*}
S_N (x) = (-1)^N \sum_{k} \bigg[{N \atop k}\bigg] (-x)^k, \tag{3.8}
\end{equation*}
where, by formula (1.18a),
\begin{equation*}
S_N^{(k)} (x) = [k]! \bigg[{N \atop k} \bigg] S_{N-k} (x), \tag{3.9}
\end{equation*}
we get
\begin{equation*}
S_N (x) = \sum_k \bigg[{N \atop k}\bigg] G_{N-k} (x -^{\hskip-.09truein \cdot}
 1)^k = \sum_{k} 
\bigg[{N \atop k}\bigg] (x -^{\hskip-.09truein \cdot}1)^{N-k} G_k, \tag{3.10}
\end{equation*}
where, by the Gauss formula (1.7), 
\begin{equation*}
G_k = S_k (1) = \left\{
\begin{array}{cl}
0, & k \ {\rm{odd}}, \\
(q^{k-1} ; q^{-2})_{\lfloor k/2 \rfloor}, & k \ {\rm{even}} .
\end{array} \right.
\tag{3.11}
\end{equation*}
Thus,
\begin{equation*}
 S_{N} (x) = \sum_k \bigg[{N \atop 2k}\bigg] (x -^{\hskip-.09truein \cdot} 1)
^{N-2k} (q^{2k-1}; q^{-2})_k . \tag{3.12}
\end{equation*}
Comparing formulae (1.12) and (3.12), we see that we must have
\begin{equation*}
\bigg[{N \atop 2k}\bigg]_q (q^{2k-1}; q^{-2})_k  =
 \left\{
\begin{array}{ll}
\bigg[{\dps{{m \atop k}}}\bigg]_{q^{2}} (q^{2m+1}; q^{-2})_k, & \ \ N = 2 m+1 
\\[1.1em]
\bigg[{\dps{{m \atop k}}}
\bigg]_{q^{2}} (q^{2m-1}; q^{-2})_k , &  \ \ N = 2m 
\end{array} \right. \tag{3.13}
\end{equation*}
and these relations can be easily verified.  Thus,
\begin{equation*}
(-1)^N \sum^N_{k=0} \bigg[{N \atop k}\bigg] (-x)^k = \sum^{\lfloor N/2 \rfloor }_{k=0} 
\bigg[{N \atop 2k}\bigg] (x -^{\hskip-.09truein \cdot} 1)^{N-2k} (q^{2k-1}; q^{-2})_k. 
\tag{3.14}
\end{equation*}
\begin{remark}{3.15}
Euler's formula (1.13) suggests that one should consider more
general family of polynomials:
\begin{equation*}
P_N (x) = \sum^{N}_{\ell = 0} \bigg[{N \atop \ell} \bigg] x^\ell q^{\alpha\ell^{2}}, 
\tag{3.16}
\end{equation*}
with $\alpha = 0$ corresponding to the Gauss case, $\alpha = 1/2$ corresponding to the Euler 
case, and $\alpha = 1$ corresponding to the Szeg${\rm{\ddot{o}}}$ case [1,7].  Applying the 
arguments used above, we find:
\begin{gather*}
{d P_N(x) \over d_q x} = [N] q^\alpha P_{N-1} (q^{2 \alpha} x), \tag{3.17}
\\
P_{N+1} (x) = q^\alpha x P_N (q^{2\alpha} x) + P_N (qx) \tag{3.18a}
\\
\qquad\qquad= q^{N + \alpha} x P_N (q^{2\alpha -1} x) + P_N (x), \tag{3.18b}
\\
P_N (x) = \sum^N_{k=0} \bigg[{N \atop k} \bigg] \rho_{N-k} (x) \theta_k, \tag{3.19}
\end{gather*}
where
\begin{equation*}
\rho_n (x) = q^{(1-2\alpha){n \choose 2}} (-q^{(2n-1)\alpha} x ; q^{-1} )_n \tag{3.20}
\end{equation*}
satisfies the same $q$-differential equation (3.17) as $P_n (x)$:
\begin{equation*}
{d \rho_n (x) \over d_q x} = [n] q^\alpha \rho_{n-1} (q^{2\alpha} x), \tag{3.21}
\end{equation*}
and $\theta_k$'s are some $x$-independent connection coefficients.  Unfortunately, I haven't been 
able to find a compact expression for the coefficients $\theta_k = \theta_k (q; \alpha)$. 
\end{remark}

\section{The geometric progressions point of view}

Formula (1.2) 
\begin{equation*}
\sum^N_{\ell=0} {N \choose \ell} (-1)^\ell = \delta^N_0 , \ \ N \in \Z_+, 
\tag{4.1}
\end{equation*} 
can be equivalently put into the following interesting form:
\begin{equation*}
\sum^\infty_{\ell =0} {t^\ell \over (1 + t)^{\ell + 1}} = 1. \tag{4.2}
\end{equation*}
(We treat all series as formal power series, and so don't have to pay attention to questions 
of convergence.  The series (4.2) converges for real $t > - 1/2$.)  Indeed, multiply equality 
(4.1) by $(-t)^N$ and then sum on all $N \in \Z_+$:
\begin{align*}
1 &= \ \sum_{N, \ell} (-t)^N {N \choose \ell} (-1)^\ell = \sum_{s, \ell} (-t)^{s+\ell} 
{s + \ell \choose \ell} (-1)^\ell = \sum_{\ell \geq 0} t^\ell \sum_{s \geq 0} 
{s + \ell \choose \ell} (-t)^s 
\\
&= \sum_{\ell \geq 0} {t^\ell \over (1 + t)^{\ell +1}} , 
\end{align*}
where we used the following version of the Newton's binomial
\begin{equation*}
{1 \over (1-t)^{N+1}} = \sum_{s \geq 0} {N + s \choose s} t^s. \tag{4.3}
\end{equation*}

We can perform similar conversion upon the formula (1.5), an Euler-type $q$-analogue of formula 
(4.1).  Multiply the equality
\begin{equation*}
\sum^N_{\ell = 0} \bigg[{N \atop \ell}\bigg] (-1)^\ell q^{{\ell \choose 2}} = \delta^N_0, \ \ 
N \in \Z_+, \tag{4.4}
\end{equation*}
by $(-t)^N$ and sum over all $N \in \Z_+$:
\begin{align*}
1 &= \sum_{N, \ell} (-t)^N \bigg[{N \atop \ell} \bigg] (-1)^\ell q^{{\ell \choose 2}} = 
\sum_{s, \ell \geq 0} (-t)^{s+\ell} \bigg[{s + \ell \atop \ell}\bigg] (-1)^\ell 
q^{{\ell \choose 2}} 
\\
&= \sum_\ell t^\ell q^{{\ell \choose 2}} \sum_s \bigg[{s + \ell \atop \ell}\bigg] (-t)^s 
\ {\rm{[by}} \ (4.6)] = \sum_{\ell \geq 0} {t^{\ell} q^{{\ell \choose 2}} \over 
(1 \dot + t)^{\ell + 1}} \ . 
\end{align*}
Thus, 
\begin{equation*}
\sum^\infty_{\ell = 0} {t^\ell q^{{\ell \choose 2}} \over (1 \dot + t ) 
^{\ell +1}} = 1; \tag{4.5}
\end{equation*}
we used in the calculation above the following Euler version of formula (4.3):
\begin{equation*}
{1 \over ( 1-^{\hskip-.09truein \cdot} t)^{N + 1}} = \sum_{s \geq 0} \bigg[{N + s \atop 
s} \bigg] t^s. \tag{4.6}
\end{equation*}

Let us now apply the same conversion device to the Gauss result (1.7):
\begin{equation*}
G_N = \sum^N_{k=0} \bigg[{N \atop k}\bigg] (-1)^k = \left\{
\begin{array}{cl}
0, & N \ {\rm{odd}}, \\
(q^{N-1} ; q^{-2})_{\lfloor N/2 \rfloor}, & N \ {\rm{even}} .
\end{array} \right.
\tag{4.7}
\end{equation*}
Multiplying by $(-t)^N$ and summing on $N$ we find:
\begin{align*}
\sum_N (-t)^N G_N &= \sum_m t^{2m} (q^{2m-1}; q^{-2})_m = 1 + \sum^\infty_{m=1} 
(1-q) ... (1 - q^{2m-1}) t^{2m}
\\
&= \sum_N (-t)^N \sum_k \bigg[{N \atop k} \bigg] (-1)^k = \sum_{s,k} (-t)^{k+s} 
\bigg[{k + s \atop k} \bigg] (-1)^k 
\\
&= \sum_k t^k \sum_s \bigg[{k+s \atop k} \bigg] (-t)^s 
= \sum_{k \geq 0} {t^k \over (1 \dot + t)^{k+1}} . 
\end{align*}
Thus,
\begin{equation*}
\sum_{k \geq 0} {t^k \over (1 \dot + t)^{k+1}} = 1 + \sum^\infty_{m=1} (1 - q) ... 
(1 - q^{2m-1}) t^{2m}. \tag{4.8}
\end{equation*}
This formula is the first from a pair found by Carlitz in [3].  The second formula in that 
pair is the case $\{r=1\}$ of the following general relation
\begin{equation*}
\sum^\infty_{\ell = 0} {(q^r t)^{\ell}q^{{\ell \choose 2}} \over
 (1 \dot + t)^{{\ell +1}} } = \sum_{N \geq 0} 
( 1 -^{\hskip-.09truein \cdot} q^r)^N (-t)^N, \tag{4.9}
\end{equation*}
which can be proven as follows:
\begin{align*}
\sum^\infty_{\ell = 0} {(q^r t)^\ell q^{{\ell \choose 2}} \over (1 \dot + t)^{\ell +1} } 
&= \sum_\ell q^{r \ell} t^\ell q^{{\ell \choose 2}} \sum_s \bigg[{\ell + s \atop \ell}\bigg] 
(-t)^s = \sum_{N \geq 0} t^N \sum^N_{\ell = 0} (-1)^{N - \ell} \bigg[{N \atop \ell}\bigg] 
q^{{\ell \choose 2}} (q^r)^\ell 
\\
&= \sum_N (-t)^N \sum^N_{\ell = 0} \bigg[{N \atop \ell}\bigg] q^{{\ell \choose 2}} 
(-q^r)^\ell \ {\rm{[by}} \ (1.3)] = \sum_N (-t)^N (1 -^{\hskip-.09truein \cdot} q^r)^N. 
\end{align*}

For $r = 0$, formula (4.9) becomes formula (4.5).  Since $r$ is arbitrary, replacing in 
formula (4.9) \ $ tq^r$ by another variable $z$, we get
\begin{equation*}
\sum^\infty_{\ell = 0} {z^\ell q^{{\ell \choose 2}} \over (1 \dot + t)^{\ell +1}} = 
\sum_{N \geq 0} (-1)^N (t -^{\hskip-.09truein \cdot} z)^N, \tag{4.10}
\end{equation*}
a $q$-analogue of the geometric progression formula 
\begin{equation*}
{1 \over 1+t } \sum^\infty_{\ell = 0} \bigg({z \over 1 + t}\bigg)^\ell
 = \sum^\infty_{N = 0} 
(z - t)^N. \tag{4.11}
\end{equation*}

\section{Gauss-like non-alternating sums}

For $x = - 1$, Newton's formula (1.1) yields
\begin{equation*}
\sum^N_{\ell =0} {N \choose \ell} = 2^N. \tag{5.1}
\end{equation*}
Similarly, the Euler binomial (1.3) for $x = - q$ provides 
\begin{equation*}
\sum^N_{\ell = 0} \bigg[{N \atop \ell}\bigg] q^{{\ell + 1 \choose 2}} = (1 \dot + q)^N. 
\tag{5.2}
\end{equation*}
If we apply to these two banalities Gauss-like ansatz, we should look at the sums of 
the form
\begin{equation*}
\sum^N_{\ell = 0} \bigg[{N \atop \ell}\bigg] (q^r)^\ell. \tag{5.3}
\end{equation*}

Not much is known about such sums, at least as far as I can tell.  (See Remark 6.12.)  
However, we shall see  
below that for $r = 1/2$, 
\begin{equation*}
\sum^N_{\ell = 0} \bigg[{N \atop \ell}\bigg] q^{\ell/2} = 
(- q^{1/2}; q^{1/2})_N. \tag{5.4}
\end{equation*}
Changing $q$ into $q^2$, this formula may be rewritten in the form
\begin{equation*}
\sigma_N = \sum^N_{\ell = 0} \bigg[{N \atop \ell}\bigg]_{q^{2}} q^\ell = (1 \dot +q)^N. 
\tag{5.5}
\end{equation*}

Let's prove it.  This formula is obviously true for $N = 0, 1.$  Using induction on $N$ 
and observing that 
\begin{equation*}
\sigma_N = \sum^N_{\ell = 0} \bigg[{N \atop \ell}\bigg]_{q^{2}} q^\ell = \sum^N_{\ell = 0} 
\bigg[{N \atop N - \ell}\bigg]_{q^{2}} q^\ell = \sum^N_{\ell =0} 
 \bigg[{N \atop \ell}\bigg]_{q^{2}} 
q^{N-\ell} , \tag{5.6}
\end{equation*}
we find:
\begin{align*}
\sigma_{N+1} = \sum_{\ell \geq 0} \bigg[{N + 1 \atop \ell}\bigg]_{q^{2}} q^\ell 
\ {\rm{[by \ (2.8b)]}} = \sum_{\ell = 0} \bigg( \bigg[{N \atop \ell}\bigg]_{q^{2}} + 
\bigg[{N\atop \ell -1}\bigg]_{q^{2}} q^{2N+2-2\ell} \bigg) q^\ell 
\\
= \sigma_N + \sum_{\ell \geq 0} \bigg[{N \atop \ell}\bigg]_{q^{2}} q^{2N+1-\ell} 
\ {\rm{[by}} \ (5.6)] = \sigma_N + q^{N+1} \sigma_N = (1+q^{N+1}) \sigma_N. \tag{5.7}
\end{align*}
Thus, 
\begin{equation*}
\sigma_{N+1} = (1 + q^{N+1}) \sigma_N, \tag{5.8}
\end{equation*} 
and since $\sigma_0 = 1$, formula (5.5) follows.

The {\it{derivation}} of formula (5.7) above suggests consideration of more general sums
\begin{equation*}
\sigma_N (\gamma) = \sum^N_{k = 0} \bigg[{N \atop k}\bigg]_{q^{2}} q^{\gamma k}. 
\tag{5.9}
\end{equation*}
Since
\begin{equation*}
\sum^N_{k=0} \bigg[{N \atop k}\bigg]_{q^{2}} q^{\gamma k} = \sum^N_{k = 0} \bigg[{ N 
\atop N - k} \bigg]_{q^{2}} q^{\gamma k} = \sum^N_{k = 0} \bigg[{N \atop k}\bigg]_{q^{2}} 
q^{\gamma (N-k)} = q^{\gamma N} \sigma_N (-\gamma), 
\end{equation*}
we find that
\begin{equation*}
\sigma_N (-\gamma) = q^{- \gamma N} \sigma_N (\gamma). \tag{5.10}
\end{equation*}
Further,
\begin{align*}
\sigma_{N+1} (\gamma) &= \sum^N_{k = 0} \bigg[{N+1 \atop k}\bigg]_{q^{2}} q^{\gamma k} 
\ {\rm{[by \ (2.8a)]}} = \sum^N_{k=0} (q^{2k} \bigg[{N \atop k}\bigg]_{q^{2}} + 
\bigg[{N \atop k-1}\bigg]_{q^{2}} ) q^{\gamma k} 
\\
&= \sigma_N (\gamma + 2) + q^\gamma \sigma _N (\gamma) , 
\end{align*}
so that
\begin{equation*}
\sigma_N (\gamma + 2) = \sigma_{N+1} (\gamma) - q^\gamma \sigma_{N} (\gamma). 
\tag{5.11}
\end{equation*}
Since we have already calculated $\sigma_N = \sigma_N (1)$ \ (5.5), formula (5.11) allows 
us to find $\sigma_N (\gamma)$ for arbitrary odd $\gamma$.

Setting
\begin{equation*}
\sigma_N (2 \ell + 1) = \sigma_N (1) \sum^\ell_{s=0} c_{\ell|s} q^{{s+1 \choose 2}} 
Q^s, \ \ Q = q^N, \ \ \ell \in \Z_+, \tag{5.12}
\end{equation*}
we can translate the recurrence relation (5.11) into the form 
\begin{equation*}
c_{\ell + 1|s} = (q^s - q^{2\ell+1} ) c_{\ell |s} + c_{\ell | s-1}, \tag{5.13}
\end{equation*}
with the understanding that
\begin{equation*}
c_{\ell | s} = 0 \ \ {\rm{unless}} \ \ 0 \leq s \leq \ell. \tag{5.14}
\end{equation*}

Since
\begin{equation*}
c_{0|0} = 1, \tag{5.15}
\end{equation*}
a little calculation shows that
\begin{gather*}
c_{\ell| 2 r} = \bigg[{\dps{\ell - r \atop r}}\bigg]_{q^{2}}
 {g_{\ell - r} \over g_r}, \tag{5.16a}
\\
c_{\ell | 2 r+1} = \bigg[{\dps{\ell - r - 1 \atop r}}\bigg]_{q^{2}}
 {g_{\ell - r} \over g_{r+1}}, \tag{5.16b}
\end{gather*}
where $g_i$'s are the Gauss products:
\begin{equation*}
g_i = {\sqcap_{t \ {\rm{odd}} \ < 2i}} (1 - q^t), \ i \in \N; \ \ g_0=1. \tag{5.17}
\end{equation*}

It's easy to verify that formulae (5.16) satisfy the recurrence relation (5.13) and the boundary 
condition (5.15).  It's interesting to observe that formula (5.16) exhibits still another form 
of 2-periodicity.

The first few $\sigma_N (2 \ell + 1)$'s are written below:
\begin{align*}
\sigma_N (3) / \sigma_N (1) &= (1-q)  + q Q, \tag{5.18a}
\\
\sigma_N (5) / \sigma_N (1) &= (1 -q) (1-q^3) + qQ (1-q^3) + q^3 Q^2, \tag{5.18b}
\\
\sigma_N (7) / \sigma_N (1) &= (1 - q) (1-q^3) (1 -q^5) + qQ (1-q^3) (1-q^5) \\
&\quad + q^3Q^2 (1-q^3) [2]_{q^{2}} + q^6 Q^3, \tag{5.18c}
\\
\sigma_N (9) / \sigma_N (1) &= (1-q) (1-q^3) (1-q^5) (1-q^7) + qQ (1-q^3) (1-q^5) (1-q^7) 
\\
&\quad {}+ q^3 Q^2 (1-q^3) (1-q^5) [3]_{q^{2}} + q^6 Q^3 (1-q^5) [2]_{q^{2}} + q^{10} Q^4. 
\tag{5.18d}
\end{align*}

Passing to the limit $N \rightarrow \infty$ and considering $|q|<1$, so that 
$Q=q^N \rightarrow 0$, we find:
\begin{equation*}
\lim_{N \to \infty} \sigma_N (2 \ell +1)/ \sigma_N (1) = (1-q) (1-q^3) ... (1-q^{2 \ell -1}), 
\ \ \ell \in \N. \tag{5.19}
\end{equation*}
Since
\begin{equation*}
\sigma_\infty (\gamma) = \lim_{N \to \infty}  \sigma_N (\gamma) = \sum_{k \geq 0} 
\bigg[{\infty \atop k}\bigg]_{q^{2}} q^{\gamma k} = 1 + \sum_{k > 0} 
{q^{\gamma k} \over (1 - q^2) ... (1 - q^{2k})} , \tag{5.20}
\end{equation*}
formula (5.19) can be rewritten as
\begin{equation*}
\sum_{k \geq 0} {q^{2\ell +1)k} \over (q^2; q^2)k} = (q; q^2)_\ell \sum_{k \geq 0} 
{q^k \over (q^2; q^2)_k} . \tag{5.21}
\end{equation*}
Now 
\begin{equation*}
(a; \rho)_\ell = (a; \rho )_\infty / (\rho^\ell a; \rho)_\infty, \tag{5.22}
\end{equation*}
so that formula (5.21) can be rewritten as
\begin{equation*}
{1 \over (q; q^2)_\infty} \sum^\infty_{k=0} {z^k \over (q^2; q^2)_k} = 
{1 \over (z; q^2)_\infty} \sum^\infty_{k=0} {q^k \over (q^2; q^2)_k} , \tag{5.23}
\end{equation*}
where we introduced
\begin{equation*}
z = q^{2\ell + 1}. \tag{5.24}
\end{equation*}
Formula (5.23) is true as it stands, for {\it{arbitrary}} $z$, because the difference 
of the LHS and the RHS of this formula is an analytic function of $z$ for $|z|<1$, vanishing 
for an infinite number of different values $z = q^{2\ell+1}, \ \ell \in \Z_+$, 
condensing to zero.
\begin{remark}{5.25} 
The {\it{alternating}} Gauss-like sums (1.9) 
\begin{equation*}
(-1)^N s_{N|r} = \sum^N_{\ell = 0} \bigg[{N \atop \ell}\bigg] (-1)^\ell (q^r)^\ell 
\tag{5.26}
\end{equation*}
have been effectively calculated in Section 1 for {\it{integer}} $r \in \Z$. \ The {\it{
non-alternating}} sums (5.3) 
\begin{equation*}
\sum^N_{\ell = 0} \bigg[{N \atop  \ell}\bigg] (q^r)^\ell \tag{5.27}
\end{equation*}
have been effectively calculated in this section for {\it{half-integers}} $r \in {1 \over 2} 
+ \Z$.  There must be some underlying reasons for this dichotomy.
\end{remark}

\section{Remarks}

\begin{remark}{6.1}
The basic philosophy of $q$-language is {\it{multiplicative}} 
discretization of classical continuous mathematics.  Interestingly enough, the formulae in 
this 
paper can be interpreted as statements in an {\it{additive}} discrete language, a certain 
$q$-analogue of the classical difference calculus.  The latter can be summarized as follows.

Let $\pmtheta = (\theta (0), \theta (1), ...)$ be a fixed sequence.  For every sequence 
$\{a_n\}$, define the $q$-difference sequences
\begin{gather*}
(\Delta^0 a)_n = a_n, \tag{6.1a}
\\
(\Delta^{k+1} a)_n = ( \Delta^k a)_{n+1} - q^{\theta(k)} (\Delta^k a)_n, \ \ k \in \Z_+. 
\tag{6.1b}
\end{gather*}
When the parameter $\pmtheta$ has the canonical form 
\begin{equation*}
\theta (k) = k, \ \ k \in \Z_+, \tag{6.2}
\end{equation*}
the sequences $\{( \Delta^ka)_n |k, n \in \Z_+ \}$ can be reconstructed from the boundary 
conditions
\begin{equation*}
b_k = (\Delta^k a)_0, \ \ \ k \in \Z_+, \tag{6.3}
\end{equation*}
by the easily verifiable formula
\begin{equation*}
(\Delta^k a)_n = \sum^n_{s=0} b_{k+n-s} \bigg[{n \atop s}\bigg] q^{ks}. \tag{6.4}
\end{equation*}

In particular, when $k = 0$ we get
\begin{equation*}
a_n = (\Delta^0 a)_n = \sum^n_{s=0} b_{n-s} \bigg[{n \atop s}\bigg] = \sum^n_{s=0}
b_s \bigg[{n \atop s}\bigg]. \tag{6.5}
\end{equation*}
Thus, evaluation of the sums (5.26) and (5.27):
\begin{equation*}
\sum^N_{\ell=0} \bigg[{N \atop \ell}\bigg] (\pm q^r)^\ell, \tag{6.6}
\end{equation*}
can be thought of as the process of reconstruction of the original sequence $\{a_N\}$ given the 
boundary $q$-difference sequence  $\{(\Delta^n a)_0 = (\pm q^r)^n\}$. 

In a superficially more general direction, say for the nonalternating case, if we  fix 
$r, \rho \in \Z_+$ and set
\begin{equation*}
b_s = \bigg[{s \atop \rho }\bigg]
 q^{\alpha (s)} \ \ \ , \ \ \ \alpha (s) = (s - \rho) (r + {1 \over 2}), \tag{6.7}
\end{equation*}
we find
\begin{align*}
a_n &=   \sum^n_{s=0} 
\bigg[{n \atop s}\bigg] b_s = \sum^n_{s=0} \bigg[{n \atop s}\bigg] 
\bigg[{s \atop \rho}\bigg] q^{\alpha (s)} = \bigg[{n \atop \rho}\bigg]
\sum^n_{s= \rho} \bigg[{n - \rho \atop s-\rho}\bigg] q^{\alpha(s) } 
\\
&= \bigg[{n \atop \rho}\bigg]
 \sum^{n-\rho}_{s=0} \bigg[{n-\rho \atop s}\bigg]q ^{s(r+{1\over 2})} = \bigg[{n \atop 
\rho}\bigg] 
\tilde \sigma_{n-\rho} 
(2r+1), \tag{6.8}
\end{align*}
where
\begin{equation*}
\tilde \sigma_N (\gamma; q) = \sigma_N (\gamma; q^{1\over 2}). \tag{6.9}
\end{equation*}
In particular, for $r=0$ and $\rho = 1$, formula (6.8) yields:
\begin{equation*}
a_n = [n] (-q^{{1\over 2}}; q^{{1\over 2}})_{n-1}. \tag{6.10}
\end{equation*}

When $q=1$, this becomes S. Rabinowitz's Crux 946 formula ([6], p. 194 )
\begin{equation*}
a_n = n \cdot 2^{n-1}, \ \ \ b_n = n, \ \ \ n \in \Z_+. \tag{6.11}
\end{equation*}
\end{remark}
\begin{remark}{6.12}
Many formulae in this paper can be found in the literature.  
The polynomials $(-1)^N S_N (-x) $ (1.11) are called by Andrews ``Rogers-Szeg${\ddot{\rm{o}}}$  
polynomials'', and many of their interesting properties are listed on pp. 49-51 in [2].  
Andrews also provides a very short proof of the Gauss  formulae (1.7), on p. 37 in [2].  
N. J. Fine has also studied these polynomials; formula (5.5) can be found on p. 29 
of his book [4], as well as on p. 49 of the Andrews book [2].
\end{remark}
\begin{remark}{6.13}
The Gauss device can be thought of as chopping off the 
naturally occurring factors $q^{{n \choose 2}}$ from the Euler $q$-analogue (1.32) of Newton's 
binomial (1.1).  In the opposite spirit, one can ask about what happens when we {\it{attach}} 
these factors to a place that is naturally missing them, another Euler's form of Newton's 
binomial,  formula (4.6):
\begin{equation*}
V_N (t) = \sum_{s \geq 0} \bigg[{N + s \atop s}\bigg] t^s q^{{s \choose 2}}. \tag{6.14}
\end{equation*}
Since these objects are no longer polynomials but are in fact infinite series, we won't 
pursue this avenue here and leave it to the reader as an exercise.  The numbers $v_N = V_N (q)$ 
can be found on p. 8 of Fine's book [4]:
\begin{gather*}
v_{2k} = {1 \over (q^2; q^2)_k} \sum_{n \geq 0} q^{{n + 1 \choose 2}} = {1 \over (q^2; 
q^2)_k} \sqcap_{n \geq 1} \bigg({1 - q^{2n} \over 1 - q^{2n-1} } \bigg), \tag{6.15a}
\\
 v_{2k+1} = {1 \over (q; q^2)_k} = {1 \over (1-q)(1-q^3)...(1 - q^{2k+1})}. \tag{6.15b}
\end{gather*}
\end{remark}

\label{lastpage}

\end{document}